\documentclass[12pt,leqno]{article}
\usepackage{color}
\usepackage{float}
\usepackage{soul}
\usepackage{amsmath}
\usepackage{amssymb}
\usepackage{theorem}
\usepackage{fancyhdr}
\usepackage{enumitem}
\usepackage{color}
\usepackage{empheq}

\usepackage[hidelinks]{hyperref}

\usepackage[margin=1in]{geometry}

\definecolor{labelkey}{rgb}{0.2,0.2,0.2}
\definecolor{refkey}{rgb}{0.2,0.2,0.2}

\newcommand{\nnn}{\ensuremath{{n\in{\mathbb N}}}}

\newcommand{\menge}[2]{\big\{{#1}~\big |~{#2}\big\}}

\newcommand{\fenv}[1]%
{\ensuremath{\,\overrightarrow{\operatorname{env}}_{#1}}}
\newcommand{\benv}[1]%
{\ensuremath{\,\overleftarrow{\operatorname{env}}_{#1}}}
\newcommand{\emp}{\ensuremath{\varnothing}}

\newcommand{\scal}[2]{\left\langle{#1},{#2}\right\rangle}

\newcommand{\RR}{\ensuremath{\mathbb R}}

\newcommand{\RP}{\ensuremath{\mathbb R}_+}

\newcommand{\argmin}{\ensuremath{\operatorname*{argmin}}}

\newcommand{\reli}{\ensuremath{\operatorname{ri}}}

\newcommand{\cone}{\ensuremath{\operatorname{cone}}}
\newcommand{\lspan}{\ensuremath{\operatorname{span}}}
\newcommand{\aff}{\ensuremath{\operatorname{aff}}}

\newcommand{\Fix}{\ensuremath{\operatorname{Fix}}}
\newcommand{\Id}{\ensuremath{\operatorname{Id}}}

\DeclareMathAlphabet{\mbbsl}{U}{bbm}{m}{sl}

\newcommand{\nc}[2]{N^{#2}_{#1}}

\newcommand{\pnX}[1]{N^{\rm prox}_{#1}} %proximal with B=X
\def\ve{\varepsilon}
\def\dd{\delta}

\def\lm{\lambda}
\def\O{\Omega}
\def\ox{{\overline x}}
\def\oth{\overline{\theta}}

\def\oy{{\overline y}}

\def\gg{\gamma}

\def\rate{\kappa}

\def\tto{\rightrightarrows}

\DeclareMathAlphabet{\mathbbsl}{U}{bbm}{m}{sl}
\def\ball{\mathbbsl{B}}

{\begin{list}{}{%
\settowidth{\labelwidth}{\textrm{#1~}}%
\setlength{\leftmargin}{\labelwidth+\labelsep}}}%requires macro calc.sty
{\end{list}}
\newtheorem{theorem}{Theorem}[section]
\newtheorem{lemma}[theorem]{Lemma}

\newtheorem{proposition}[theorem]{Proposition}
\newtheorem{definition}[theorem]{Definition}
\theoremstyle{plain}{\theorembodyfont{\rmfamily}
}
\theoremstyle{plain}{\theorembodyfont{\rmfamily}
}
\theoremstyle{plain}{\theorembodyfont{\rmfamily}
}
\theoremstyle{plain}{\theorembodyfont{\rmfamily}
\newtheorem{example}[theorem]{Example}}
\newtheorem{fact}[theorem]{Fact}
\theoremstyle{plain}{\theorembodyfont{\rmfamily}
\newtheorem{remark}[theorem]{Remark}}
\def\proof{\noindent{\it Proof}. \ignorespaces}
\def\endproof{\ensuremath{\hfill \quad \square}}

\definecolor{labelkey}{rgb}{0,0.08,0.45}
\definecolor{refkey}{rgb}{0,0.6,0.0}
\definecolor{Brown}{rgb}{0.45,0.0,0.05}
\definecolor{lime}{rgb}{0.00,0.8,0.0}
\definecolor{lblue}{rgb}{0.5,0.5,0.99}

\definecolor{myblue}{rgb}{.9, .95, 1}
  \newcommand*\mybluebox[1]{%
    \colorbox{myblue}{\hspace{1em}#1\hspace{1em}}}

\def\doi{DOI}
\newcommand{\papercited}[6]
{{#1}, #2, \emph{#3} {\bf #4} (#5), #6.}
\newcommand{\others}[2]
{{#1}, #2}
\newcommand{\bookcited}[4]
{{#1}, \emph{#2}, #3 (#4).}
\allowdisplaybreaks % or locally if problems {\allowdisplaybreaks
\begin{document}

\title{Linear Convergence of the Douglas--Rachford Method for Two Closed Sets}
\author{Hung M.~Phan\thanks{Department of Mathematical Sciences, University of Massachusetts Lowell, Lowell, MA 01854, USA.\newline
\indent\!\;\ \
Email: \texttt{hung\_phan@uml.edu}}}

\date{February 19, 2015}

\maketitle \thispagestyle{fancy}

%\vskip 8mm

\begin{abstract} \noindent
In this paper, we investigate the Douglas--Rachford method for two closed (possibly nonconvex) sets in Euclidean spaces. We show that under certain regularity conditions, the Douglas--Rachford method converges locally with $R$-linear rate. In convex settings, we prove that the linear convergence is global. Our study recovers recent results on the same topic.
\end{abstract}


{\small \noindent {\bfseries 2010 Mathematics Subject
Classification:} 
{Primary
49M27, %% Decomposition methods
65K10; %% Optimization and variational techniques
Secondary 
47H09, %% Contraction-type mappings, nonexpansive mappings, $A$-proper mappings, etc.
49J52, %% Nonsmooth analysis
%49M20, %% Methods of relaxation type
49M37, %% Methods of nonlinear programming type
65K05, %% Mathematical programming methods
90C26. %% Nonconvex programming, global optimization
%% Nonlinear programming
%% 94A08. %% Image processing (compression, reconstruction, etc.)
}}

%%\renewcommand{\baselinestretch}{1.1}
%%%%%%%%%%%%%%%%%%%%%%%%%%%%%%%%%%%%%%%%%%%%%%%%%%
\section{Introduction}
Throughout the paper, 
\begin{empheq}[box=\mybluebox]{equation}\label{e:stand}
\text{$X$ is a finite dimensional Euclidean space.}
\end{empheq}
Let $A$ and $B$ be two closed subsets of $X$. The basic feasibility problem is to 
\begin{equation}\label{e:feas}
\mbox{find a point in} \ A\cap B.
\end{equation}
This problem has long been considered very important in the natural sciences and engineering. The reference \cite{BB93siam} is often considered a classic survey of methods for solving \eqref{e:feas}. Among them, the Douglas--Rachford method \cite{dr56ams} has attracted increasing attention, mainly because of its good performance.

To describe this method, we first recall that the {\em distance function} to a closed subset $A$ of $X$ is $d_A(x):=\inf_{a\in A}\|x-a\|$; the {\em projector} $P_A$ and {\em reflector} $R_A$ are the set-valued mappings defined respectively by
\begin{align}
&P_A(x):=\argmin_{a\in A}\|x-a\|=\menge{a\in A}{\|x-a\|=d_A(x)}\\
\mbox{and}\quad
&R_A x:=(2P_A-\Id)(x)=\menge{2a-x}{a\in P_A x}.
\end{align}
When $P_A(x)=\{a\}$ is a singleton, we simply write $a=P_Ax$.

The Douglas-Rachford operator $T:X\tto X$ for two sets $A$ and $B$ is then defined by
\begin{equation}\label{e:DRgs}
T(x):=\menge{P_{B}(2a-x)+x-a}{a\in P_{A}x}=\tfrac{1}{2}(\Id+R_{B} R_{A})(x).
\end{equation}
Clearly when $A$ and $B$ are convex, $T$ is single-valued.

Next, let $x_0$ be a starting point, the {\em Douglas--Rachford method (DR)} for two sets $A$ and $B$ generates sequences $(x_n)_\nnn$ such that
\begin{equation}
\forall\nnn:\quad x_{n+1}\in T(x_n).
\end{equation}
Each such sequence $(x_n)_\nnn$ is called a DR sequence with starting point $x_0$.  

If the DR sequence converges to a fixed point $\ox$, then there exists an element of $P_A\ox$ that is a solution of \eqref{e:feas}. Thus, DR can be used to solve \eqref{e:feas}.

It is well-known that for convex subsets of a Hilbert spaces, the DR sequence converges weakly to a fixed point, see, e.g., \cite{LM79,svai11}. In nonconvex settings, the theoretical justification is still incomplete. In fact, the paper \cite{BorSim11} proves the local convergence for DR between a sphere and an affine line; while \cite{AB12} proves the global convergence of DR for a line and a circle in $\RR^2$. It is worth mentioning that, DR has also been showed to be useful in sparse affine feasibility problem \cite{{DeZh09mit},{HLN13b}}.

In this paper, we focus on the $R$-linear convergence of DR sequences. We say that a sequence $(x_n)_\nnn$ converges to a point $\bar{x}$ with \emph{$R$-linear rate} $\rate\in[0,1)$ if there exists $C\geq 0$ such that
\begin{equation}
\forall\nnn:\quad
\|x_n-\bar{x}\|\leq C\rate^n.
\end{equation}

Recently, Hesse and Luke \cite{HL13siopt} have obtained an interesting result about the local $R$-linear convergence for DR of two sets in nonconvex settings. In particular, the authors proved that: ``{\em if $A$ is an affine subspace and $B$ is a superregular set (see Definition~\ref{d:eps-rg}), and the system $\{A,B\}$ is strongly regular (see~\eqref{e:st-rg}), then the DR sequence converges locally to the intersection $A\cap B$ with $R$-linear rate}" (see \cite[Theorem~3.18]{HL13siopt}).

We will complement the above statement with several new results:
\begin{enumerate}[label=(R\arabic*)]
\item\label{rs:i0} If $A$ and $B$ are two superregular sets, and the system $\{A,B\}$ is strongly regular (see~\eqref{e:st-rg}), then the DR sequence converges locally with {\em $R$-linear rate} to the intersection $A\cap B$ (see Theorem~\ref{t:m1}).

\item\label{rs:i} If $A$ and $B$ are two superregular sets, and the system $\{A,B\}$ is affine-hull regular (see~\eqref{e:af-rg}), then the DR sequence converges locally with {\em $R$-linear rate} to a fixed point $\ox$ of $T$ and that $P_A\ox\equiv P_B\ox\in A\cap B$, which solves the feasibility problem \eqref{e:feas} (see Theorem~\ref{t:m2}).

\item\label{rs:ii} If $A$ and $B$ are two convex sets such that $\reli A\cap\reli B\neq\emp$, then for every starting point, the DR sequence converges with {\em $R$-linear rate} to a fixed point of $T$ and that $P_A\ox\equiv P_B\ox\in A\cap B$, which solves the feasibility problem \eqref{e:feas} (see Theorem~\ref{t:m3}).
\end{enumerate}

Several comments are in order:
\begin{enumerate}
\item \ref{rs:i0} is more general than \cite[Theorem~3.18]{HL13siopt} (see Example~\ref{ex:2cir}).

\item \ref{rs:i} is more general than \ref{rs:i0} (see~Example~\ref{ex:thm2>thm1}).

\item In \ref{rs:i0}, the limit point $\ox$ of the DR sequence is indeed a solution of the feasibility problem \eqref{e:feas}.
\item In \ref{rs:i} and \ref{rs:ii}, although the limit point $\ox$ may not be a solution of \eqref{e:feas}, the ``shadow" $P_A\ox\equiv P_B\ox\in A\cap B$ is actually a solution. Notice that, we implicitly assume that both projectors $P_A$ and $P_B$ are computable. Therefore, as long as the limit point $\ox$ is obtained, the shadows $P_A\ox$ and $P_B\ox$ are computable.

\item In \ref{rs:i} and \ref{rs:ii}, the limit point $\ox$ is a fixed point of the DR operator and surprisingly, $P_A\ox\equiv P_B\ox$ is a singleton. However, for a general DR fixed point $\ox$, $P_A\ox$ does not necessarily coincide with $P_B\ox$ (see Example~\ref{ex:PA=PB}) and that neither of them is necessarily singleton.

\item Although the theory of DR for two convex sets is well-known, the {\em rate of convergence} for this case has only been observed partially before: the case of two affine subspaces \cite{BBNPW14,HLN13b}, the case of one convex set and one affine subspace \cite{HL13siopt}. We would like to mention that \ref{rs:ii} is the {\em first} to address the {\em $R$-linear convergence} of DR for two closed convex sets. This result is also re-established in \cite[Section~4]{BNP15} within the context of averaged nonexpansive operators.
\end{enumerate}

The paper is organized as follows. Section~\ref{s:pre} contains preliminary results. In Section~\ref{s:afre}, we present an {\em affine reduction property} of DR operator. Section~\ref{s:M} then presents the main results of the paper. Finally, concluding remarks are given in Section~\ref{s:7}.

\medskip
\noindent{\bf Notation:} The notation used in this paper is standard and follows \cite{{BC2011},{Boris1},{Rock98}}. The smallest affine subspace containing a set $\O$ is denoted by $\aff\O$. The relative interior of $\O$, $\reli\O$, is the interior of $\O$ relative to $\aff\O$. We write $\Phi:X\tto X$ if $\Phi$ is a multi-valued mapping, and we denote $\Fix \Phi:=\menge{x\in X}{x\in \Phi x}$ the set of fixed points of $\Phi$. For $z\in X$ and $\dd\in\RP$, \ $\ball_\dd(z):=\menge{x\in\RR}{\|x-z\|\leq \dd}$ is the closed ball centered at $z$ with radius $\dd$. We also write $\ball:=\ball_1(0)$ and $\ball_\dd:=\ball_\dd(0)$ for brevity.

%%%%%%%%%%%%%%%%%%%%%%%%%%%%%%%%%%%%%%%%%%
%%%%%%%%%%%%%%%%%%%%%%%%%%%%%%%%%%%%%%%%%%
\section{Preliminary results}
\label{s:pre}

In this section, we recall some preliminary concepts and results.

\subsection{Normal cones and regularity of systems of sets}

Recall that the {\em proximal normal cone} (see, e.g., \cite[Example~6.16]{Rock98} or \cite[eq.~(2.80)]{Boris1}) to a set $A$ at a point $x$ is defined by $\pnX{A}(x):=\cone(P^{-1}_A x-x)$. The {\em Mordukhovich (or limiting) normal cone} (see, e.g., \cite[Theorem~1.6]{Boris1}) is given by
\begin{equation}
N_A(x):=\menge{u\in X}{\exists \{x_n\}_\nnn\subset A,\ u_n\in\pnX{A}(x_n): \ x_n\to x, u_n\to u}.
\end{equation}
Let $L$ be an affine subspace containing $A$, then the (limiting) normal cone of $A$ restricted to $L$ is given by
\begin{equation}\label{e:rnc}
\nc{A}{L}(x):=N_A(x)\cap(L-x).
\end{equation}
Clearly, $\nc{A}{L}(x)\subseteq\nc{A}{}(x)$. Indeed, the concept of restricted normal cone was developed in \cite{BLPW13a,BLPW13b,BPW14} in more general settings.

%%%%%%%%%%%%%%%%%%%%%%%%%%%%%%%%%%%%%
\begin{definition}[regularity of set systems]
Let $A$ and $B$ be two subsets of $X$ and let $w\in A\cap B$. We say that the system $\{A,B\}$ is {\em strongly regular} at $w$ if
\begin{equation}\label{e:st-rg}
\nc{A}{}(w)\cap(-\nc{B}{}(w))=\{0\}.
\end{equation}
Let $L:=\aff(A\cup B)$, we say that the system $\{A,B\}$ is {\em affine-hull regular} at $w$ if
\begin{equation}\label{e:af-rg}
\nc{A}{L}(w)\cap(-\nc{B}{L}(w))=\{0\}.
\end{equation}
\end{definition}

\begin{remark}
\label{r:strg}
Clearly, strong regularity implies affine-hull regularity. However, the reverse is not always true. For example, consider two distinct lines $A$ and $B$ in $\RR^3$ that intersect at only one point $w$. Then $\{A,B\}$ is affine-hull regular at $w$, but is not strongly regular at $w$. Condition \eqref{e:st-rg} is called the normal qualification condition \cite[Definition~3.2]{Boris1}, while \eqref{e:af-rg} is called the qualification condition for systems of sets \cite[Definition~6.6]{BLPW13a}.
\end{remark}

In order to quantify the rate of convergence, we need the following quantity for two closed cones $N_1$ and $N_2$ of $X$ (cf., the CQ-number \cite{BLPW13a,BLPW13b})
\begin{equation}
\begin{aligned}
\overline{\theta}(N_1,N_2)&:=\sup\menge{\scal{u}{v}}{u\in N_1\cap\ball,v\in-N_2\cap\ball}\\
&=\max\menge{\scal{u}{v}}{u\in N_1\cap\ball,v\in-N_2\cap\ball},
\end{aligned}
\end{equation}
which is related to strong/affine-hull regularity as showed below (a simple proof is included since it uses only the definitions).

\begin{lemma}\label{l:0915a}
Let $A$ and $B$ be two closed sets and let $w\in A\cap B$. Then the following hold
\begin{enumerate}
\item\label{l:0915ai} $\{A,B\}$ is strongly regular at $w$ if and only if \ $\overline{\theta}(N_A(w),N_B(w))<1$.
\item\label{l:0915aii} $\{A,B\}$ is affine-hull regular at $w$ if and only if \ $\overline{\theta}\big(\nc{A}{L}(w),\nc{B}{L}(w)\big)<1$, where $L:=\aff(A\cup B)$.
\end{enumerate}
\end{lemma}
\proof \ref{l:0915ai}: Notice that $\overline{\theta}(N_A(w),N_B(w))\in[0,1]$, and we have
$\overline{\theta}(N_A(w),N_B(w))=1$ $\Leftrightarrow$ $\exists u\in N_A(w)\cap\ball$, $\exists v\in-N_B(w)\cap\ball$ such that $\scal{u}{v}=1$
$\Leftrightarrow$ $\exists u\in N_A(w)\cap\ball$, $\exists v\in-N_B(w)\cap\ball$ such that $u=v$ and $\|u\|=\|v\|=1$
$\Leftrightarrow$ $N_A(w)\cap(-N_B(w))\neq\{0\}$ $\Leftrightarrow$ $\{A,B\}$ is not strongly regular. This implies \ref{l:0915ai}.

\ref{l:0915aii}: We may assume $w=0$, thus $L=\aff(A\cup B)=\lspan(A\cup B)$ is a subspace. Then restrict the consideration to the subspace $L$.
\endproof

Indeed, strong regularity (or more general, uniform regularity) can be characterized by different quantities \cite{kru-th-13}. Thorough discussions can be found in \cite{kru06, kru-th-13} and the references therein.

%%%%%%%%%%%%%%%%%%%%%%%%%%%%%%%%%
Finally, we finish this section by recalling linear regularity property of set systems \cite{BBL99}, which plays an important role in convex and variational analysis. Its connection to strong regularity is also stated below.

%%%%%%%%%%%%%%%%%%%%%%%%%%%%%%%%%
\begin{definition}[linear regularity of set systems {\cite{BBL99}}]
\label{d:linrg}
We say that the system of sets $\{A,B\}$ is {\em $\mu$-linear regular} on $S$ if
\begin{equation}
\forall x\in S:\quad d_{A\cap B}(x)\leq 
\mu\max\{ d_{A}(x),d_{B}(x)\}.
\end{equation}
We also say that $\mu$ is a linear regularity modulus.
\end{definition}

%%%%%%%%%%%%%%%%%%%%%%%%%%%%%%%%%
\begin{fact}
{\rm (\cite[p.~2410]{HL13siopt}, and \cite[Theorem~1]{kru06})}
\label{f:kru} 
Suppose that $\{A,B\}$ is strongly regular at $w\in A\cap B$. Then there exist $\dd>0$ and $\mu\geq 1$ such that $\{A,B\}$ is $\mu$-linear regular on $\ball_\dd(w)$, i.e.,
\begin{equation}
\forall x\in \ball_\dd(w):\quad d_{A\cap B}(x)\leq 
\mu\max\{ d_{A}(x),d_{B}(x)\}.
\end{equation}
\end{fact}

%%%%%%%%%%%%%%%%%%%%%%%%%%%%%%%%%%%%%%%%%%%%
\subsection{Regularity of sets and quasi firm nonexpansiveness}

It is well-known that if $A$ and $B$ are convex, then their DR operator $T$ is firmly nonexpansive (see, e.g., \cite[Proposition~4.21]{BC2011}). In particular, the following holds
\begin{equation}\label{e:frm}
\forall x\in X,\forall \ox\in\Fix T:\quad
\|x-Tx\|^2+\|Tx-\ox\|^2\leq\|x-\ox\|^2.
\end{equation}
If $A$ and $B$ are not convex, however, \eqref{e:frm} is not necessarily true.

Nevertheless, we will show that an analogous estimation for $T$ holds locally under certain regularity assumptions on the sets $A$ and $B$ (see Proposition~\ref{p:0914b}). First, we recall some technical definitions.

\begin{definition}[superregularity]
\label{d:eps-rg}
Let $A$ be a closed subset of $X$. We say that $A$ is {\em $(\ve,\dd)$-regular} at $w$ if $\ve\geq0$, $\dd>0$ and
\begin{equation}\label{e:eps-rg}
\left.\begin{aligned}
&x,z\in A\cap \ball_\dd(w),\\
& u\in\pnX{A}(x)
\end{aligned}\right\}
\quad\Rightarrow\quad
\scal{u}{z-x}\leq\ve\|u\|\cdot\|z-x\|.
\end{equation}
We say that $A$ is {\em superregular} at $w$ if for every $\ve>0$, there exists $\dd>0$ such that $A$ is $(\ve,\dd)$-regular at $w$.
\end{definition}
Superregularity was first introduced in \cite[Definition~4.3]{LLM09}. It is somewhat between Clarke regularity and amenability or prox-regularity. This concept is further generalized in \cite[Definition~8.1]{BLPW13a} and \cite[Definition~2.9]{HL13siopt}. Importantly, all convex sets are superregular. More discussion and examples can be found in \cite{BLPW13a,HL13siopt,LLM09}.

%%%%%%%%%%%%%%%%%%%%%%%%%%%%%%%%%%%%%%
\begin{definition}[$\gg$-quasi firm expansiveness]
\label{d:qfne}
A mapping $\Phi:X\tto X$ is said to be {\em $(\O,\gg)$-quasi firmly nonexpansive} on $U$ if for all $x\in U$, $x_+\in \Phi x$, and $\ox\in P_\O x$, we have
\begin{equation}
\|x-x_+\|^2+\|x_+-\ox\|^2\leq \gg \|x-\ox\|^2.
\end{equation}
\end{definition}

One would notice that when $\gg=1$, quasi firm expansiveness is a variant of the quasi nonexpansiveness in \cite[Definition~4.1]{BC2011}. In addition, the latter is a restriction of Fej\'er monotonicity \cite[Chapter~5]{BC2011}. Besides, the $\gg$-quasi firm expansiveness is a simplification of the $(\O,\ve)$-firm nonexpansiveness \cite[Definition~2.3(ii)]{HL13siopt}: a mapping $\Phi:X\tto X$ is said to be $(\O,\ve)$-firmly nonexpansive on $U$ if
\begin{equation}
\begin{aligned}
\forall x\in U,\ &\forall x_+\in \Phi x,\ \forall \ox\in\O,\ \forall\ox_+\in\Phi\ox:\\
&\|x_+-\ox_+\|^2+\|(x-x_+)-(\ox-\ox_+)\|^2\leq (1+\ve) \|x-\ox\|^2.
\end{aligned}
\end{equation}

Although the following argument is simplified from \cite{HL13siopt}, details are included for the readers' convenience.

%%%%%%%%%%%%%%%%%%%%%%%%%%%%%%%%%%%%%%%%%%%%%%
\begin{fact}{\rm(\cite[Theorem~2.14]{HL13siopt})}
\label{f:1227c} 
Let $A$ be a closed subset of $X$ and suppose that $A$ is $(\ve,\dd)$-regular at $w\in A$. Let $U:=\menge{x\in X}{P_{A}x\subset\ball_\dd(w)}$. Then
\begin{enumerate}
\item\label{f:1227ci}
For all $x\in U$, $a\in P_{A}x$, $\ox\in A\cap\ball_\dd(w)$:\quad $\|a-\ox\|\leq(1+\ve)\|x-\ox\|$.
\item\label{f:1227ciii}
For all $x\in U$, $u\in R_Ax$, $\ox\in A\cap\ball_\dd(w)$:\quad $\|u-\ox\|\leq(1+2\ve)\|x-\ox\|$.
\end{enumerate}
\end{fact}

%%%%%%%%%%%%%%%%%%%%%%%%%%%%%%%%%%%%%%%%%%%%%%

\begin{lemma}
\label{l:1229a} Let $A$ and $B$ be two closed sets.
Let $\dd>0$, $\ve\in[0,\tfrac{1}{4})$. Assume $A$ is $(\ve,2\dd)$-regular at $w\in A\cap B$. Then the following hold:
\begin{enumerate}
\item\label{l:1229ai} 
$\forall x\in \ball_\dd(w):$\quad
$P_{A}x\subset\ball_{2\dd}(w)$.

\item\label{l:1229aii}
$\forall x\in \ball_\dd(w):$\quad
$P_{B}R_{A}x\subset\ball_{3\dd}(w)$.
\end{enumerate}
\end{lemma}
%----------
\proof
Take an $x\in \ball_\dd(w)$. \ref{l:1229ai}: For every $a\in P_{A}x$, we have $\|a-w\|\leq\|a-x\|+\|x-w\|\leq 2\|x-w\|\leq 2\dd$. So, $P_{A}x\subseteq \ball_{2\dd}(w)$. 

\ref{l:1229aii}: Applying Fact~\ref{f:1227c}\ref{f:1227ciii} to $A$ being $(\ve,2\dd)$-regular at $w$, we have
\begin{equation}
\forall u\in R_{A}x:\quad
\|u-w\|\leq(1+2\ve)\|x-w\|.
\end{equation}
For every $b\in P_{B}R_{A}x$,
let $u\in R_{A}x$ such that $b\in P_{B}u$. So
\begin{equation}
\|b-w\|\leq\|b-u\|+\|u-w\|\leq 2\|u-w\|
\leq 2(1+2\ve)\dd\leq 3\dd,
\end{equation}
which implies \ref{l:1229aii}.
\endproof

%%%%%%%%%%%%%%%%%%%%%%%%%%%%%%%%%%%%%%%%%%%%%%%%%%%%%%%%%
\begin{proposition}
\label{p:0914b}
Let $A$ and $B$ be two closed sets and let $T$ be the DR operator \eqref{e:DRgs}. Let $\dd>0$, $\ve_1,\ve_2\in[0,\tfrac{1}{4})$.
Assume further that $A$ and $B$ are $(\ve_1,2\dd)$- and $(\ve_2,3\dd)$-regular at $w\in A\cap B$, respectively. Define also
\begin{equation}
\O:=A\cap B%%\cap\ball_{2\dd}(w)
\quad\mbox{and}\quad
\gg:=\tfrac{1+(1+2\ve_1)^2(1+2\ve_2)^2}{2}.
\end{equation}
Then the following hold
\begin{enumerate}
\item\label{p:0914bi}
For all $x\in \ball_\dd(w)$, $x_+\in Tx$, and $\ox\in \O\cap\ball_{2\dd}(w)$, we have
\begin{equation}
\|x-x_+\|^2+\|x_+ -\ox\|^2\leq \gg\|x-\ox\|^2.
\end{equation}

\item\label{p:0914bii} For all $x\in \ball_\dd(w)$, $x_+\in Tx$, and $\ox\in P_\O x$, we have
\begin{equation}
\|x-x_+\|^2+\|x_+ -\ox\|^2\leq \gg\|x-\ox\|^2,
\end{equation}
i.e., $T$ is $(\O,\gg)$-quasi firmly nonexpansive on $\ball_\dd(w)$ (see Definition~\ref{d:qfne}).
\end{enumerate}
\end{proposition}
%----------
\proof Take any $x\in \ball_\dd(w)$ and $x_+\in Tx$.

\ref{p:0914bi}: We can find $u\in R_{A}x$, $b\in P_{B}u$, and $v=2b-u$ such that $x_+=\tfrac{1}{2}(x+v)$. Now take any $\ox\in \O\cap\ball_{2\dd}(w)$. First, applying Fact~\ref{f:1227c}\ref{f:1227ciii} to $A$ being $(\ve_1,2\dd)$-regular, we have
\begin{equation}\label{e:p0914b1}
\|u-\ox\|\leq(1+2\ve_1)\|x-\ox\|.
\end{equation}
By Lemma~\ref{l:1229a}\ref{l:1229aii}, we have $b\in\ball_{3\dd}(w)$. Then applying Fact~\ref{f:1227c}\ref{f:1227ciii} to $B$ being $(\ve_2,3\dd)$-regular at $w$, we have
\begin{equation}
\|v-\ox\|\leq (1+2\ve_2)\|u-\ox\|.
\end{equation}
Combining with \eqref{e:p0914b1},
\begin{equation}\label{e:p0914b2}
\|v-\ox\|\leq(1+2\ve_1)(1+2\ve_2)\|x-\ox\|.
\end{equation}
Next, using the parallelogram law, we obtain
\begin{equation}
2\|x-\ox\|^2+2\|v-\ox\|^2=\|x-v\|^2+\|2(x_+-\ox)\|^2
=4\|x-x_+\|^2+4\|x_+-\ox\|^2
\end{equation}
Employing \eqref{e:p0914b2}, we have
\begin{equation}\label{e:p0914b3}
\|x-x_+\|^2+\|x_+-\ox\|^2
=\tfrac{1}{2}(\|x-\ox\|^2+\|v-\ox\|^2)\leq\Big(\tfrac{1+(1+2\ve_1)^2(1+2\ve_2)^2}{2}\Big)\|x-\ox\|^2.
\end{equation}

\ref{p:0914bii}: Take any $\ox\in P_\O x$, we have
\begin{equation}
\|\ox-w\|\leq\|\ox-x\|+\|x-w\|\leq 2\|x-w\|\leq 2\dd.
\end{equation}
So $\ox\in \O\cap\ball_{2\dd}(w)$. Thus, the conclusion follows from \ref{p:0914bi}.
\endproof

%%%%%%%%%%%%%%%%%%%%%%%%%%%%%%%%%%%%%%%%%%%
Finally, we include an $R$-linear convergence result of Fej\'{e}r monotonicity type with an elementary proof for completeness.

%%%%%%%%%%%%%%%%%%%%%%%%%%%%%%%%%%%%%%
\begin{proposition}\label{p:0104a}
Let $\Phi:X\tto X$ be an operator, let $\O$ be a closed subset of $X$ and $w\in\O$. Assume that there are $\dd>0$ and $\rate\in[0,1)$ such that
\begin{equation}\label{e:p0104a1}
\forall x\in \ball_\dd(w),\ \forall x_+\in\Phi x,\ \forall \ox\in P_\O x:\quad\|x_+-\ox\|\leq\rate \|x-\ox\|.
\end{equation}
Let $x_0$ and let \ $x_{n+1}\in\Phi x_n$ \ for all $\nnn$. Then if $x_0$ is sufficiently close to $w$, the sequence $(x_n)_\nnn$ converges $R$-linearly to a point $\ox$. In particular,
\begin{equation}\label{e:p0104a3}
\|x_n-\ox\|\leq \tfrac{\|x_0-w\|(1+\rate)}{1-\rate}\rate^n
\quad\text{and}\quad
\ox\in\O\cap\ball_\dd(w).
\end{equation}
\end{proposition}
%----------
\proof
Assume that $x_0$ satisfies $M:=\|x_0-w\|\leq \frac{\dd(1-\rate)}{2}$, we will show that the conclusion holds. Indeed, let $\ox_n\in P_\O x_n$ for all $\nnn$. We first prove by induction that
\begin{equation}\label{e:p0104a2}
\forall \nnn:\quad 
\|x_{n+1}-\ox_{n}\|\leq\rate\|x_{n}-\ox_{n}\|\leq M\rate^{n+1}.
\end{equation}
It is easy to check that \eqref{e:p0104a2} holds for $n=0$. Now, suppose \eqref{e:p0104a2} holds for $0,\ldots,n-1$, we will show it also holds for $n$. Indeed,
\begin{equation}
\begin{aligned}
\|x_n-w\|
&%\leq \|x_0-w\|+\sum_{i=0}^{n-1}\|x_{i+1}-x_{i}\|
\leq \|x_0-w\|+\sum_{i=0}^{n-1}\big(\|x_{i+1}-\ox_{i}\| +\|x_{i}-\ox_{i}\|\big)\\
&\leq M+\sum_{i=0}^{n-1}\big(M\rate^{i+1}+M\rate^{i}\big)
\leq 2M\sum_{i=0}^n\rate^i
\leq 2M\tfrac{1}{1-\rate}\leq \dd.
\end{aligned}
\end{equation}
So by \eqref{e:p0104a1} applied to $x_n\in\ball_\dd(w)$, $x_{n+1}\in\Phi x_n$, and $\ox_n\in P_\O x_n$, we have
\begin{equation}
\|x_{n+1}-\ox_n\|\leq \rate \|x_n-\ox_n\|=\rate d_\O(x_n)
\leq \rate \|x_n-\ox_{n-1}\|\leq M\rate^{n+1}.
\end{equation}
So \eqref{e:p0104a2} holds for $n$, thus, it holds for all $\nnn$ by the principle of mathematical induction.

Next, using \eqref{e:p0104a2} for $0<n<m$,
\begin{equation}\label{e:p0104a4}
\begin{aligned}
\|x_n-x_m\|
&\leq \sum_{i=n}^{m-1}\big(\|x_{i+1}-\ox_{i}\| +\|x_{i}-\ox_{i}\|\big)
\leq \sum_{i=n}^{m-1}\big(M\rate^{i+1}+M\rate^{i}\big)
\leq \tfrac{M(1+\rate)}{1-\rate}\rate^n.
\end{aligned}
\end{equation}
So $(x_n)_\nnn$ is a Cauchy sequence, thus, it converges to a limit point $\ox\in\ball_\dd(w)$. Since $d_\O(x_n)=\|x_n-\ox_n\|\to 0$, we conclude that $\ox\in\O$. Finally, by letting $m\to\infty$ in \eqref{e:p0104a4}, we obtain the estimate in \eqref{e:p0104a3}.
\endproof

%%%%%%%%%%%%%%%%%%%%%%%%%%%%%%%%%%%%%
%%%%%%%%%%%%%%%%%%%%%%%%%%%%%%%%%%%%%
\section{Affine reduction}
\label{s:afre}

The main result of this section is Theorem~\ref{t:1212a}, in which we prove an interesting property of DR for two closed sets. First, we recall some properties of affine subspaces. Recall that two affine subspaces $L_1$ and $L_2$ are {\em parallel} if there exists $z\in X$ such that $L_1=z+L_2$. The proof of the following lemma is elementary, thus, omitted.
%%%%%%%%%%%%%%%%%%%%%%%%%%%%%%%%%%%%%
\begin{lemma}\label{l:1207a}
Let $L$ be an affine subspace of $X$ and let $\O\subseteq L$ be closed. Then the following hold
\begin{enumerate}
\item\label{l:1207ai} $R_\O P_L=P_LR_\O$.
\item\label{l:1207ai2} For every $x\in X$ and $z\in R_\O x$, we have: \ $x-P_Lx=P_Lz-z$.
\item\label{l:1207aii} {\rm(\cite[Corollary~3.20]{BC2011})}
$P_L$ is an affine operator, i.e., for all $x,y\in X$, $\lm\in\RR$,
\begin{equation}
P_L\big((1-\lm)x+\lm y\big)=(1-\lm)P_Lx+\lm P_Ly.
\end{equation}
%\item\label{l:0121ci} Let $L_0$ be a subspace of $X$. For every $z\in X$, we have
%\begin{equation}
%P_{(z+L_0)}(\cdot)=P_{L_0}(\cdot)+P_{L_0^\bot}z.
%\end{equation}
%\item\label{l:0121cii} Let $L_1$ and $L_2$ be two parallel affine subspaces. Then \ $P_{L_1}P_{L_2}=P_{L_1}$.
\end{enumerate}
\end{lemma}

%%%%%%%%%%%%%%%%%%%%%%%%%%%%%%%%%%%%%%%%%%%%%%%%
\begin{lemma}\label{l:0913a}
Let $A$ and $B$ be two closed sets with $L:=\aff(A\cup B)$, and let $T$ be the DR operator \eqref{e:DRgs}. Let $x_0\in X$ and $x_1\in Tx_0$. Define also $y_0=P_Lx_0$ and $y_1=P_Lx_1$. Then
\begin{equation}
y_1\in Ty_0\quad
\text{and}\quad x_1-y_1=x_0-y_0.
\end{equation}
\end{lemma}
\proof
Since $x_1\in Tx_0$, we find $z\in R_Ax_0$ and $w\in R_Bz$ such that $x_1=\tfrac{1}{2}(x_0+w)$. Employing Lemma~\ref{l:1207a}\ref{l:1207ai}, we have
\begin{equation}
P_Lw\in P_L R_B R_Ax_0=R_B P_L R_Ax_0= R_B R_A P_L x_0=R_BR_A y_0.
\end{equation}
Now since $P_L$ is an affine operator (see Lemma~\ref{l:1207a}\ref{l:1207aii}), we have
\begin{equation}
y_1=P_Lx_1=\tfrac{1}{2}P_Lw+\tfrac{1}{2}P_Lx_0\in\tfrac{1}{2}(R_BR_A y_0+y_0)=Ty_0.
\end{equation}
This proves the first part.

To prove the second part, we employ Lemma~\ref{l:1207a}\ref{l:1207ai2} twice to obtain
\begin{equation}
w-P_Lw=P_Lz-z=x_0-P_Lx_0=x_0-y_0.
\end{equation}
So, $x_1-y_1=\frac{1}{2}(x_0+w)-P_L\big(\frac{1}{2}(x_0+w)\big)
=\frac{1}{2}(x_0-P_Lx_0)+\frac{1}{2}(w-P_Lw)=x_0-y_0$.
\endproof

%%%%%%%%%%%%%%%%%%%%%%%%%%%%%%%%%%%%%%%%%%%%%%%%

\begin{theorem}\label{t:1212a}
Let $A$ and $B$ be two closed sets with $L:=\aff(A\cup B)$, and let $T$ be the DR operator \eqref{e:DRgs} Let $(x_n)_\nnn$ be a DR sequence with a starting point $x_0\in X$. Let $y_n:=P_L x_n$ for $\nnn$. Then
\begin{enumerate}
\item\label{t:1212ai} $\forall\nnn)$ \ $x_n-y_n=x_0-y_0$.
\item\label{t:1212aii} $(y_n)_\nnn$ is a DR sequence generated by $T$ with starting point $y_0$.
\item\label{t:1212aiii} If $(y_n)_\nnn$ converges to $\oy \in A\cap B$, then $(x_n)_\nnn$ converges to $\ox$ such that $y_n-\oy=x_n-\ox$. In particular, the rates of convergence are identical. Moreover,
\begin{equation}\label{e:t1212aiii}
%\quad\mbox{and}\quad
P_A\ox=P_B\ox=P_L\ox=\ox-(x_0-y_0)=\oy,
\end{equation}
which is a singleton.
\end{enumerate}
\end{theorem}
%----------
\proof \ref{t:1212ai}\&\ref{t:1212aii}: apply Lemma~\ref{l:0913a} inductively.

\ref{t:1212aiii}: Suppose $(y_n)_\nnn$ converges to $\oy\in A\cap B$, then \ref{t:1212ai} implies $(x_n)_\nnn$ converges to
\begin{equation}
\ox=\oy+(x_0-y_0)=\oy+(x_n-y_n).
\end{equation}
This also implies $y_n-\oy=x_n-\ox$. Thus, the rates of convergence are identical.

Next, since $P_L(\cdot)$ is continuous, we have $\oy=P_L \ox$. On the other hand, $\oy\in A\subseteq L$, so $\oy=P_A\ox$ is a singleton. Similarly, $\oy=P_B\ox$. Thus, we obtain \eqref{e:t1212aiii}.
\endproof
%%%%%%%%%%%%%%%%%%%%%%%%%%%%%%%%%%%%%%%%%%

\begin{remark}\label{r:1213a}
Theorem~\ref{t:1212a} is a surprising result in the sense that, in order to study DR for two sets $A$ and $B$, it suffices to study DR with starting points in $L=\aff(A\cup B)$. Indeed, let $(x_n)_\nnn$ be a DR sequence and let $y_n=P_Lx_n$. Theorem~\ref{t:1212a}\ref{t:1212ai} implies that the behaviors of the sequence $(y_n)_\nnn$ will imply the behaviors of the sequence $(x_n)_\nnn$ and vice vesa.
\end{remark}

%%%%%%%%%%%%%%%%%%%%%%%%%%%%%%%%%%%%%%%
%%%%%%%%%%%%%%%%%%%%%%%%%%%%%%%%%%%%%%%
\section{Main results}
\label{s:M}
This section contains the main results of the paper which are divided into three parts.

%%%%%%%%%%%%%%%%%%%%%%%%%%%%%%%%%%%%%%%
\subsection{DR under strong regularity}
\label{ss:strg}
In this part, we establish the $R$-linear convergence of DR for two sets $A$ and $B$ locally around $w\in A\cap B$ under two assumptions
\begin{enumerate}
\item $\{A,B\}$ is strongly regular at $w$ (see~\eqref{e:st-rg});\quad and
\item $A$ and $B$ are superregular at $w$ (see Definition~\ref{d:eps-rg}).
\end{enumerate}
The result of this section (Theorem~\ref{t:m1} below) is an improvement of \cite[Theorem~3.18]{HL13siopt} where the authors proved the local $R$-linear convergence assuming one set is superregular and the other is an affine subspace.

\begin{lemma}\label{l:0914a}
Let $A$ and $B$ be two closed sets. Assume that $\{A,B\}$ is strongly regular at $w\in A\cap B$, or equivalently (see Lemma~\ref{l:0915a}),
\begin{equation}
\oth:=\oth(N_A(w),N_B(w))<1.
\end{equation}
Then for every $\theta\in(\overline{\theta},1)$, there exist $\dd>0$ such that
\begin{equation}
\left.\begin{aligned}
a\in A\cap \ball_\dd(w),\ b\in B\cap \ball_\dd(w),\\
u\in \pnX{A}(a),\ v\in\pnX{B}(b)
\end{aligned}\right\}
\ \Longrightarrow\
\scal{u}{v}\geq -\theta\|u\|.\|v\|.
\end{equation}
\end{lemma}
\proof Suppose on the contrary that there exist sequences $a_n\in A$, $b_n\in B$, $a_n\to w$, $b_n\to w$, $u_n\in\pnX{A}(a_n)$, $v_n\in\pnX{B}$ such that
\begin{equation}
\scal{u_n}{v_n}<-\theta\|u_n\|.\|v_n\|.
\end{equation}
By dividing by $\|u_n\|.\|v_n\|$, we can assume that $u_n$, $v_n$ are unit vectors. So let $u$ and $v$ be accumulation (unit) vectors of $u_n$ and $v_n$ respectively, we have $u\in N_A(w)$, $v\in N_B(w)$, and $\scal{u}{v}\leq-\theta$. So by the definition of $\oth$,
\begin{equation}
\oth\geq\scal{u}{-v}\geq \theta>\oth,
\end{equation}
which is a contradiction.
\endproof

The following lemma provides the main ingredient.
\begin{lemma}
\label{l:m1}
Let $A$ and $B$ be two closed sets, and let $T$ be the DR operator \eqref{e:DRgs}. Assume further that $A$ is superregular at $w$, and that $\{A,B\}$ is strongly regular at $w$, or equivalently,
\begin{equation}\label{e:lm1a}
\oth:=\oth(N_A(w),N_B(w))<1.
\end{equation}
Then for any $\theta\in(\oth,1)$, there exist $\dd>0$ and $\lambda\in[0,1)$ such that for all $x\in \ball_\dd(w)$ and $x_+\in Tx$, one has
\begin{equation}\label{e:lm1h}
\|x-x_+\|
\geq
\lm
\max\big\{d_{A\cap B}(2a-x),
\tfrac{1}{\sqrt{5}}d_{A\cap B}(x)\big\}
\geq
\tfrac{\lm}{\sqrt{5}}
d_{A\cap B}(x),
\end{equation}
where $a\in P_Ax$ is such that $x_+\in P_B(2a-x)+x-a$.
\end{lemma}
%---------
\proof 
Take an arbitrary $\oth<\theta<1$, $\ve\in[0,\tfrac{1}{4})$ and define $\O:=A\cap B$. By the superregularity of $A$, Lemma~\ref{l:0914a}, Fact~\ref{f:kru}, we find $\dd>0$ and $\mu\geq 1$ such that
\begin{enumerate}
\item $A$ is $(\ve,2\dd)$-regular at $w$;
\item\label{e:lm1e}
$\left\{\begin{aligned}
a\in A\cap \ball_{2\dd}(w),\ b\in B\cap \ball_{3\dd}(w),\\
\zeta_1\in \pnX{A}(a),\ \zeta_2\in\pnX{B}(b)
\end{aligned}\right.
\ \Longrightarrow\
\scal{\zeta_1}{\zeta_2}\geq-\theta\|\zeta_1\|.\|\zeta_2\|$;
\quad and
\item\label{e:lm1g} $\forall x\in\ball_{2\dd}(w):\quad
d_\O(x)\leq\mu\max\{d_A(x),d_B(x)\}$.
\end{enumerate}

Take any $x\in\ball_\dd(w)$ and $x_+\in Tx$. Take also any $a\in P_{A}x$, $u=2a-x$, and $b\in P_{B}(u)$ such that
\begin{equation}
x_+=b+x-a\in P_{B}(u)+x-a=P_{B}(2a-x)+x-a.
\end{equation}
Since $A$ is $(\ve,2\dd)$-regular at $w$, Lemma~\ref{l:1229a} implies $a\in\ball_{2\dd}(w)$ and $b\in\ball_{3\dd}(w)$. Also, Fact~\ref{f:1227c}\ref{f:1227ciii} implies
\begin{equation}\label{e:lm1i}
\|u-w\|\leq(1+2\ve)\|x-w\|\leq 2\dd.
\end{equation}
Next, we have
\begin{equation}
a-u=x-a\in\pnX{A}(a)
\quad\mbox{and}\quad
u-b\in \pnX{B}(b).
\end{equation}
So \ref{e:lm1e} implies
\begin{equation}\label{e:lm1d}
\scal{a-u}{u-b}\geq -\theta\|a-u\|.\|u-b\|.
\end{equation}
We now have
\begin{subequations}\label{e:lm1c}
\begin{align}
\|x-x_+\|^2&=\|a-b\|^2
=\|a-u\|^2+\|u-b\|^2+2\scal{a-u}{u-b}\\
&\geq \|a-u\|^2+\|u-b\|^2-2\theta\|a-u\|.\|u-b\|\\
&\geq \|a-u\|^2+\|u-b\|^2-\theta(\|a-u\|^2+\|u-b\|^2)\\
&=(1-\theta)(\|a-u\|^2+\|u-b\|^2).
\end{align}
\end{subequations}
On the one hand, \eqref{e:lm1i} and \ref{e:lm1g} imply
\begin{equation}\label{e:lm1b}
\|a-u\|^2+\|u-b\|^2\geq\max\{d_{B}^2(u),d_{A}^2(u)\}
\geq\tfrac{1}{\mu^2}d_{\O}^2(u).
\end{equation}
On the other hand, triangle and Cauchy-Schwarz inequalities imply
\begin{subequations}
\begin{align}
d_{B}^2(x)&\leq(\|x-u\|+d_{B}(u))^2
=(2\|a-u\|+\|u-b\|)^2\\
&\leq (2^2+1^2)(\|a-u\|^2+\|u-b\|^2).
\end{align}
\end{subequations}
Using \ref{e:lm1g} again, we have
\begin{subequations}\label{e:lm1f}
\begin{align}
d_{\O}^2(x)
&\leq \mu^2\max\{d_{A}^2(x),d_{B}^2(x)\}\\
&\leq \mu^2\max\{\|a-u\|^2,5(\|a-u\|^2+\|u-b\|^2)\}\\
&\leq 5\mu^2(\|a-u\|^2+\|u-b\|^2).
\end{align}
\end{subequations}
Combining \eqref{e:lm1c}, \eqref{e:lm1b}, and \eqref{e:lm1f}, we have
\begin{equation}
\|x-x_+\|^2\geq
\Big(\tfrac{1-\theta}{\mu^2}\Big)\max
\big\{d_{\O}^2(u),
\tfrac{1}{5}d_{\O}^2(x)\big\},
\end{equation}
which yields \eqref{e:lm1h} with $\lm:=\tfrac{\sqrt{1-\theta}}{\mu}\in[0,1)$.
\endproof

%%%%%%%%%%%%%%%%%%%%%%%%%%%%%%%%%%%%%%%

%%%%%%%%%%%%%%%%%%%%%%%%%%%%%%%%%%%%%%%
\begin{theorem}[main result 1]
\label{t:m1}
Let $A$ and $B$ be two closed sets, and let $T$ be the DR operator \eqref{e:DRgs}. Assume further that $A$ and $B$ are superregular at $w\in A\cap B$, and that $\{A,B\}$ is strongly regular at $w$, or equivalently,
\begin{equation}
\oth:=\oth(N_A(w),N_B(w))<1.
\end{equation}
Then if the starting point $x_0$ is sufficiently closed to $w$, the generated DR sequence $(x_n)_\nnn$ converges to a point $\ox\in A\cap B$ with $R$-linear rate.
\end{theorem}
\proof
Define $\O:=A\cap B$. First, applying Lemma~\ref{l:m1}, there exist $\dd>0$ and $\lm\in[0,1)$ such that for all $x\in\ball_\dd(w)$ and $x_+\in Tx$, we have
\begin{equation}\label{e:tm1a}
\|x-x_+\|\geq\tfrac{\lm}{\sqrt{5}}d_{\O}(x).
\end{equation}
Take $\ve_1,\ve_2\in[0,\tfrac{1}{4})$ small and shrink $\dd$ if necessary, we assume that $A$ is $(\ve_1,2\dd)$-regular at $w$ and that $B$ is $(\ve_2,3\dd)$-regular at $w$. 

Now take any $x\in \ball_\dd(w)$, $x_+\in Tx$, and $\ox\in P_\O x$, Proposition~\ref{p:0914b}\ref{p:0914bii} implies 
\begin{equation}
\|x-x_+\|^2+\|x_+ -\ox\|^2\leq \gg\|x-\ox\|^2
\quad\text{where}\quad
\gg:=\tfrac{1+(1+2\ve_1)^2(1+2\ve_2)^2}{2}.
\end{equation}
Combining with \eqref{e:tm1a}, we have
\begin{equation}
\begin{aligned}
\|x_+-\ox\|^2
&\leq \gg \|x-\ox\|^2 -\|x-x_+\|^2\\
&\leq \gg \|x-\ox\|^2 -\lm^2\|x-\ox\|^2\\
&=(\gg-\lm^2)\|x-\ox\|^2=:\rate^2\|x-\ox\|^2.
\end{aligned}
\end{equation}
Note that $\rate^2=\tfrac{1+(1+2\ve_1)^2(1+2\ve_2)^2}{2}-\lm^2\in[0,1)$ if $\ve_1,\ve_2$ was chosen small enough.

This assures assumption \eqref{e:p0104a1} in Proposition~\ref{p:0104a} holds. Thus, the conclusion now follows from Proposition~\ref{p:0104a}. 
\endproof

\begin{remark}[rate of convergence]
\label{r:rt}
From the proofs of Proposition~\ref{p:0104a}, Lemma~\ref{l:m1}, and Theorem~\ref{t:m1}, we derive a formula for the $R$-linear rate $\rate$ as follows:

Suppose that there are $\dd,\ve_1,\ve_2>0$, $\theta\in[0,1)$, and $\mu\geq 1$ such that
\begin{enumerate}
\item $A$ is $(\ve_1,2\dd)$-regular at $w$;
\item $B$ is $(\ve_2,3\dd)$-regular at $w$;
\item
$\left\{\begin{aligned}
a\in A\cap \ball_{2\dd}(w),\ b\in B\cap \ball_{3\dd}(w),\\
\zeta_1\in \pnX{A}(a),\ \zeta_2\in\pnX{B}(b)
\end{aligned}\right.
\ \Longrightarrow\
\scal{\zeta_1}{\zeta_2}\geq-\theta\|\zeta_1\|.\|\zeta_2\|$;
\item 
$\forall x\in\ball_{2\dd}(w):\quad
d_{A\cap B}(x)\leq\mu\max\{d_A(x),d_B(x)\}$;\quad and
\item\label{r:rt-v} $\rate^2:=\tfrac{1+(1+2\ve_1)^2(1+2\ve_2)^2}{2}-\tfrac{1-\theta}{5\mu^2}\in[0,1)$.
\end{enumerate}
Then for any starting point $x_0$ such that $\|x_0-w\|\leq\tfrac{\dd(1-\rate)}{2}$, the generated DR sequence $(x_n)_\nnn$ converges to some point $\ox$ with $R$-linear rate $\rate$, more specifically,
\begin{equation}
\|x_n-\ox\|\leq \tfrac{\|x_0-w\|(1+\rate)}{1-\rate}\rate^n
\quad\text{and}\quad
\ox\in A\cap B\cap\ball_\dd(w).
\end{equation}
\end{remark}

Theorem~\ref{t:m1} is more general than \cite[Theorem~3.18]{HL13siopt}. However, in the case that $A$ is an affine subspace, we obtain a {\em smaller bound} for the $R$-linear rate. Details are given in the following result (notice that the rate $\rate$ in Theorem~\ref{t:HL318}\ref{t:HL318-iv} is smaller than $\rate$ in Remark~\ref{r:rt}\ref{r:rt-v}).

%%%%%%%%%%%%%%%%%%%%%%%%%%%%%%%%%%%%%%%%%%%%%%%%%%%%%%%
\begin{theorem}{\rm(\cite[Theorem~3.18]{HL13siopt})}
\label{t:HL318}
Let $A$ be an affine subspace, let $B$ be a closed set, and let $T$ be the DR operator \eqref{e:DRgs}. Suppose that there are $\dd,\ve,>0$, $\theta\in[0,1)$, and $\mu\geq 1$ such that
\begin{enumerate}
\item $B$ is $(\ve_2,3\dd)$-regular at $w$;
\item
$\left\{\begin{aligned}
a\in A\cap \ball_{2\dd}(w),\ b\in B\cap \ball_{3\dd}(w),\\
\zeta_1\in \pnX{A}(a),\ \zeta_2\in\pnX{B}(b)
\end{aligned}\right.
\ \Longrightarrow\
\scal{\zeta_1}{\zeta_2}\geq-\theta\|\zeta_1\|.\|\zeta_2\|$;
\item 
$\forall x\in\ball_{2\dd}(w):\quad
d_{A\cap B}(x)\leq\mu\max\{d_A(x),d_B(x)\}$;\quad and
\item\label{t:HL318-iv} $\rate^2:=\tfrac{1+(1+2\ve_1)^2(1+2\ve_2)^2}{2}-\tfrac{1-\theta}{\mu^2}\in[0,1)$.
\end{enumerate}
Then for any starting point $x_0$ such that $\|x_0-w\|\leq\tfrac{\dd(1-\rate)}{2}$, the generated DR sequence $(x_n)_\nnn$ converges to some point $\ox$ with $R$-linear rate $\rate$.
\end{theorem}
%----------
\proof
Define $\O:=A\cap B$. First, applying Lemma~\ref{l:m1}, there exist $\dd>0$ and $\lm\in[0,1)$ such that for all $x\in\ball_\dd(w)$ and $x_+\in Tx$, we have
\begin{equation}\label{e:tHL318a}
\|x-x_+\|\geq\lm
\max\big\{d_{\O}(2a-x),
\tfrac{1}{\sqrt{5}}d_{\O}(x)\big\}.
\end{equation}
Notice that since $A$ is an affine subspace and $\O\subset A$, we have $d_{\O}(2a-x)=d_{\O}(x)$. So \eqref{e:tHL318a} implies
\begin{equation}
\|x-x_+\|\geq\lm d_{\O}(x).
\end{equation}
The rest of the proof is analogous to that of Theorem~\ref{t:m1}.
\endproof

Finally, we provide a simple example in $\RR^2$ showing that Theorem~\ref{t:m1} is applicable while the results in \cite{HL13siopt} are not.

\begin{example}[DR for two transversal circles in $\RR^2$]
\label{ex:2cir}
Let $A$ and $B$ be two circles in $\RR^2$ that intersect transversally at two distinct points). Then the results in \cite{HL13siopt} are not applicable. On the other hand, since circles are superregular sets, Theorem~\ref{t:m1} is applicable and yields $R$-linear convergence for DR sequences locally around each intersection point.
\end{example}

%%%%%%%%%%%%%%%%%%%%%%%%%%%%%%%%%%%%%%%%%%%%%%%%%%%%
%%%%%%%%%%%%%%%%%%%%%%%%%%%%%%%%%%%%%%%%%%%%%%%%%%%%%%%%%%%%%%
\subsection{DR under affine-hull regularity}
\label{s:afrg}

In this part, we establish the $R$-linear convergence of DR for two sets $A$ and $B$ locally around $w\in A\cap B$ under two assumptions
\begin{enumerate}
\item $\{A,B\}$ is affine-hull regular at $w$ (see~\eqref{e:st-rg});\quad and
\item $\{A,B\}$ is superregular at $w$ (see~\eqref{e:af-rg}).
\end{enumerate}
Using Remark~\ref{r:1213a} and Theorem~\ref{t:1212a}, our strategy is to rely on the behavior of DR on the affine-hull $L=\aff(A\cup B)$.

\begin{theorem}[main result 2]
\label{t:m2}
Let $A$ and $B$ be two closed sets with $L:=\aff(A\cup B)$, and let $T$ be the DR operator \eqref{e:DRgs}. Assume further that $A$ and $B$ are superregular at $w\in A\cap B$, and that $\{A,B\}$ is affine-hull regular at $w$, i.e.,
\begin{equation}
\oth:=\oth(\nc{A}{L}(w),\nc{B}{L}(w))<1.
\end{equation}
Then if the shadow $P_Lx_0$ is sufficiently closed to $w$, the generated DR sequence $(x_n)_\nnn$ converges to a point $\ox\in\Fix T$ with $R$-linear rate. Moreover, 
\begin{equation}\label{e:shadows}
P_A\ox\equiv P_B\ox=\ox-(x_0-P_Lx_0)\in A\cap B,
\end{equation}
i.e., $P_A\ox\equiv P_B\ox$ \ solves the feasibility problem \eqref{e:feas}.
\end{theorem}
%----------
\proof 
Define $y_n:=P_Lx_n$ for $\nnn$. Then $(y_n)$ is a DR sequence with starting point $y_0=P_Lx_0$.

Now consider two set $A$, $B$ within the affine subspace $L$. Clearly, $\{A,B\}$ is affine-hull regular implies that it is strongly regular within the space $L$. Thus, Theorem~\ref{t:m1} yields that $(y_n)_\nnn$ converges to $\oy\in A\cap B$ with $R$-linear rate. In turn, Theorem~\ref{t:1212a}\ref{t:1212aiii} implies that $x_n$ converges to some $\ox$ also with $R$-linear rate, and that
\begin{equation}
P_A\ox\equiv P_B\ox\equiv P_L\ox=\oy\in A\cap B.
\end{equation}
Finally, this last relation implies $\ox\in\Fix T$.
\endproof

%%%%%%%%%%%%%%%%%%%%%%%%%%%%%%%%%%%%%
\begin{remark}
Theorem~\ref{t:m2} proves that, the region of convergence is actually {\em larger} than a ball around $w$. In fact, {\em the region of convergence is the cylinder generated by some ball around $w$ and $(L-w)^\bot$, the orthogonal complement of $L-w$}.
\end{remark}

\begin{example}\label{ex:thm2>thm1}
Theorem~\ref{t:m2} is indeed more general than Theorem~\ref{t:m1}. In fact, consider two distinct lines $A$ and $B$ in $\RR^3$ that intersect at only one point. Then one can check that Theorem~\ref{t:m2} is applicable while Theorem~\ref{t:m1} is not (recall Remark~\ref{r:strg}).
\end{example}

%%%%%%%%%%%%%%%%%%%%%%%%%%%%%%%%%%%%%%%
\subsection{DR for two convex sets}
\label{ss:cvx}

We now study the case that both sets $A$ and $B$ are convex. Because of convexity, we claim that all of the assumptions required for $R$-linear convergence will be fulfilled using only the {\em standard qualification condition} of convex analysis
\begin{equation}
\reli A\cap\reli B\neq\emp.
\end{equation}
Next, we will verify this claim.

\begin{fact}{\rm (\cite[Theorem~3.13]{BLPW13a})}
\label{f:0105b} Let $A$ and $B$ be two closed convex sets. The following are equivalent:
\begin{enumerate}
\item $\reli A\cap\reli B\neq\emp$.
\item $\{A,B\}$ is affine-hull regular at $w$ for all \ $w\in A\cap B$.
\item $\{A,B\}$ is affine-hull regular at $w$ for some \ $w\in A\cap B$.
\end{enumerate}
\end{fact}

\begin{fact}{\rm(\cite[Proposition~4.6.1]{Bthes})}
\label{f:0110a}
Let $A$ and $B$ be closed convex subsets of $X$ such that
\begin{equation}
\reli A\cap\reli B\neq\emp.
\end{equation}
Then for every bounded set $S$, there exists $\mu\geq 1$ such that $\{A,B\}$ is $\mu$-linear regular on $S$.
\end{fact}

%%%%%%%%%%%%%%%%%%%%%%%%%%%%%%
\begin{lemma}\label{l:0105a}
Let $A$ and $B$ be two closed convex sets such that $\reli A\cap\reli B\neq\emp$. Then
\begin{equation}
\big(\aff(A\cup B)\big)\cap\Fix T=A\cap B.
\end{equation}
\end{lemma}
%----------
\proof Denote $L:=\aff(A\cup B)$. ``$\supseteq$": clear. ``$\subseteq$": Let $x\in L\cap\Fix T$. Since $P_{A}$, $P_{B}$ are single-valued, let $a:=P_{A} x$, $b:=P_{B}(2a-x)$. So we have
\begin{equation}\label{e:l0105a1}
x-a\in\nc{A}{L}(a)
\quad\mbox{and}\quad
(2a-x)-b\in\nc{B}{L}(b).
\end{equation}
Since $x$ is a fixed point of $T$, $x=Tx=x+b-a$. So $a=b\in A\cap B$. Thus, it follows from \eqref{e:l0105a1} that \ $x-a\in\nc{A}{L}(a)$ \ and \ $a-x\in\nc{B}{L}(a)$. Employing Fact~\ref{f:0105b}, we have $x-a\in\nc{A}{L}(a)\cap(-\nc{B}{L}(a))=\{0\}$. Hence, $x=a=b\in A\cap B$.
\endproof

%%%%%%%%%%%%%%%%%%%%%%%%%%%%%%%%%%%%%%%%%%%%%%%%%%
\begin{fact}{\rm(\cite[Theorem~25.6]{BC2011})}
\label{f:0105c} Let $A$ and $B$ be two closed convex sets, and let $T$ be the DR operator \eqref{e:DRgs}. Let $(x_n)_\nnn$ be the DR sequence with an arbitrary starting point $x_0\in X$. Then $(x_n)_\nnn$ converges to some point in \ $\Fix T$.
\end{fact}

We are then ready to present the main result of this section.
%%%%%%%%%%%%%%%%%%%%%%%%%%%%%%%%%%%%
\begin{theorem}[two convex sets]
\label{t:m3}
Let $A$ and $B$ be two closed convex sets with $L:=\aff(A\cup B)$, and let $T$ be the DR operator \eqref{e:DRgs}. Assume also that
\begin{equation}
\reli A\cap\reli B\neq\emp.
\end{equation}
Then for every starting point $x_0\in X$, the DR sequence $(x_n)_\nnn$ converges to a point $\ox\in \Fix T$ with $R$-linear rate and that
\begin{equation}
P_A\ox\equiv P_B\ox=\ox-(x_0-P_Lx_0)\in A\cap B.
\end{equation} 
Furthermore, the shadows $(P_Ax_n)_\nnn$ and $(P_Bx_n)_\nnn$ also converge to $P_A\ox\equiv P_B\ox\in A\cap B$ with $R$-linear rate.
\end{theorem}
%----------
\proof
Let $y_n=P_Lx_n$ for $\nnn$, then Theorem~\ref{t:1212a} implies $(y_n)_\nnn$ is also a DR sequence and that $x_n-y_n=x_0-y_0$. Then Fact~\ref{f:0105c} implies that both $(x_n)_\nnn$ and $(y_n)_\nnn$ converge to the fixed points $\ox$ and $\oy$ respectively. So $\ox-\oy=x_0-y_0$. We also have $\oy\in L\cap\Fix T=A\cap B$ by Lemma~\ref{l:0105a}.

Now, since $\{A,B\}$ is affine-hull regular at $\oy$,  Theorem~\ref{t:m2} implies that $(x_n)_\nnn$ and $(y_n)_\nnn$ converges $R$-linearly and that $P_A\ox=P_B\ox=\ox-(x_0-y_0)\in A\cap B$. Finally, $(P_Ax_n)_\nnn$ also converges $R$-linearly due to the nonexpansiveness of $P_A$, similarly for $(P_Bx_n)_\nnn$.
\endproof

%%%%%%%%%%%%%%%%%%%%%%%%%%%%%%%%%%%%%
%%%%%%%%%%%%%%%%%%%%%%%%%%%%%%%%%%%%%
\section{Remarks and examples}
\label{s:7}

We conclude the paper by some remarks.

\begin{example}[two convex sets]
In $\RR^2$, consider two (convex) strips
\begin{equation}
A=\menge{(x_1,x_2)}{0\leq x_2\leq 1}
\quad\mbox{and}\quad
B=\menge{(x_1,x_2)}{0\leq x_2-x_1\leq 1}.
\end{equation}
One can check that $\{A,B\}$ is strongly regular at every point in the intersection $A\cap B$.

Notice that the results in \cite{HL13siopt} are not applicable because neither $A$ nor $B$ is an affine subspace. Theorem~\ref{t:m3}, on the other hand, does apply and yield convergence for the DR with even global $R$-linear rate.

We now use Remark~\ref{r:rt} to compute the rate: first, notice that $\{A,B\}$ is (globally) linearly regular with modulus $\mu=\tfrac{1}{\sin\tfrac{\pi}{8}}=\tfrac{2}{\sqrt{2-\sqrt{2}}}$.Next, we set $\ve_1=\ve_2=0$, $\dd=+\infty$, and $\theta=\tfrac{\sqrt{2}}{2}$. Then, the rate is
\begin{equation}
\rate^2=\frac{1+(1+2\ve_1)^2(1+2\ve_2)^2}{2}
-\frac{1-\theta}{5\mu^2}
=1-\frac{1-\tfrac{\sqrt{2}}{2}}{5\cdot\tfrac{4}{2-\sqrt{2}}}=\frac{17+2\sqrt{2}}{20}.
\end{equation}
Thus, $\rate=\sqrt{\tfrac{17+2\sqrt{2}}{20}}$. Despite the conjecture that the {\em actual} rate could be smaller, our obtained rate $\rate$ is the {\em only} rate known so far!
\end{example}

%%%%%%%%%%%%%%%%%%%%%%%%%%%%%%%%%%%%%
\begin{remark}[the sequence of interest]
In Theorem~\ref{t:m3}, the statement on the shadow sequences $(P_Ax_n)_\nnn$ and $(P_Bx_n)_\nnn$ only holds in convex settings. In nonconvex settings, the mentioned statement is not always true. To that end, the behavior of $(x_n)$ is indeed important, see also \cite[p.~2398]{HL13siopt}.

In fact, in convex settings, the sequence of interest is $(P_Ax_n)_\nnn$, rather than $(x_n)_\nnn$ itself. The reason is that the distance from shadow $P_A x_n$ to $P_A\ox$ can be {\em very small} after very few steps even though the true iteration $x_n$ is far away from the limit $\ox$. This yields the ``rippling'' behavior of the method, see \cite[Figure~1]{BBNPW14} and also \cite[Figures~4 and 6]{BKroad}. This can provide an edge over other methods, for example, {\em von-Neumann method of alternating projections}, see \cite[Section~9]{BBNPW14}.
\end{remark}

%%%%%%%%%%%%%%%%%%%%%%%%%%%%%%%%%%%%%
\begin{example}[$P_A\ox\equiv P_B\ox$ may fail without affine-hull regularity]
\label{ex:PA=PB}
In Theorems~\ref{t:m2} and \ref{t:m3}, the conclusion $P_A\ox\equiv P_B\ox$ may fail if the affine-hull regularity is violated. For example, in $\RR^2$, consider
\begin{equation}
A=\RR\times\RP,
B=\menge{(x_1,x_2)}{x_2\leq -x_1^2},
\quad\mbox{and}\quad
w=(0,0).
\end{equation}
Then $L=\aff(A\cup B)=\RR^2$. Clearly, $\nc{A}{L}(w)\cap(-\nc{B}{L}(w))=\{0\}\times\RR_-$, which means the affine-hull regularity does not hold. Now take $x_0=(0,-1)$, then $x_0\in\Fix T$, hence the DR sequence is $(x_n)\equiv x_0$. However, $P_A x_0=(0,0)$ is different from $P_B x_0=x_0$.
\end{example}

%%%%%%%%%%%%%%%%%%%%%%%%%%%%%%%%%%%%%%%
%%%%%%%%%%%%%%%%%%%%%%%%%%%%%%%%%%%%%%%
\section*{Acknowledgment}
The author was partially supported by an internal grant of University of Massachusetts Lowell. This research was initiated during the author's stay at University of British Columbia (Kelowna, Canada). The author is grateful to Heinz Bauschke and Xianfu Wang (UBC Kelowna, Canada) for their support. The author also thanks the editors and the referees for their constructive comments.

%\end{document}

%%%%%%%%%%%%%%%%%%%%%%%%%%%%%%%%%%%%%%%%%%%%%%%%%%%

\end{document}